\documentclass[titlepage,12pt]{article}

\usepackage{amssymb}
\usepackage{latexsym}

\newcommand{\ben}{\begin{enumerate}}
\newcommand{\een}{\end{enumerate}}
\newcommand{\ble}{\begin{lem}}
\newcommand{\ele}{\end{lem}}
\newcommand{\bth}{\begin{thm}}
\renewcommand{\eth}{\end{thm}}
\newcommand{\bpr}{\begin{prop}}
\newcommand{\epr}{\end{prop}}
\newcommand{\bco}{\begin{cor}}
\newcommand{\eco}{\end{cor}}
\newcommand{\bcon}{\begin{conj}}
\newcommand{\econ}{\end{conj}}
\newcommand{\bde}{\begin{defn}}
\newcommand{\ede}{\end{defn}}
\newcommand{\bex}{\begin{exa}}
\newcommand{\eex}{\end{exa}}
\newcommand{\barr}{\begin{array}}
\newcommand{\earr}{\end{array}}
\newcommand{\btab}{\begin{tabular}}
\newcommand{\etab}{\end{tabular}}
\newcommand{\beq}{\begin{equation}}
\newcommand{\eeq}{\end{equation}}
\newcommand{\bea}{\begin{eqnarray*}}
\newcommand{\eea}{\end{eqnarray*}}
\newcommand{\bce}{\begin{center}}
\newcommand{\ece}{\end{center}}
\newcommand{\bpi}{\begin{picture}}
\newcommand{\epi}{\end{picture}}
\newcommand{\bfi}{\begin{figure} \begin{center}}
\newcommand{\efi}{\end{center} \end{figure}}

\newcommand{\bsl}{\begin{slide}{}}
\newcommand{\esl}{\end{slide}}

\newcommand{\pf}{{\bf Proof.}}
\newcommand{\qed}{\rule{1ex}{1ex}}
\newcommand{\Qed}{\rule{1ex}{1ex} \medskip}

\newcommand{\hs}[1]{\hspace{#1}}

\newcommand{\fl}[1]{\lfloor #1 \rfloor}

\newcommand{\cef}[2]{\left\lceil\frac{#1}{#2}\right\rceil}
\newcommand{\gau}[2]{\left[ #1 \atop #2 \right]}

\newcommand{\ze}{\zeta}

\newcommand{\bn}{{\bf n}}

\newcommand{\fS}{{\mathfrak S}}

\newcommand{\Des}{\mathop{\rm Des}}

\newcommand{\inv}{\mathop {\rm inv}}

\newcommand{\tang}{\mathop{\rm tang}}

\newcommand{\pams}{Proc. Amer. Math. Soc.}

\newtheorem{thm}{Theorem}[section]
\newtheorem{prop}[thm]{Proposition}
\newtheorem{cor}[thm]{Corollary}
\newtheorem{lem}[thm]{Lemma}
\newtheorem{conj}[thm]{Conjecture}
\newtheorem{exa}[thm]{Example}
\newtheorem{defn}[thm]{Definition}

\begin{document}
\bibliographystyle{plain}
\pagestyle{empty}
\title{Arithmetic Properties of  Generalized Euler Numbers}
\author{ Bruce E. Sagan\\
Department of Mathematics \\ 
Michigan State University\\
East Lansing, MI 48824-1027\\[.1in]
and\\[.1in]
Ping Zhang \\
Department of Mathematics and Statistics \\ 
Western Michigan University\\ 
Kalamazoo, MI 49008-5152}
 
\date{\today \\[1in]
        \begin{flushleft}
        \end{flushleft}
       }
\maketitle

\begin{flushleft} Proposed running head: \end{flushleft}
     \begin{center}
Generalized Euler Numbers
     \end{center}

Send proofs to:
\begin{center}
Bruce E. Sagan \\ Department of Mathematics \\Michigan State
University \\ East Lansing, MI 48824-1027\\[5pt]
Tel.: 517-355-8329\\
FAX: 517-336-1562\\
Email: sagan@math.msu.edu
\end{center}

         \begin{abstract}
\quad \  The generalized Euler number $E_{n|k}$ counts
the number of permutations of $\{1,2,\ldots,n\}$ which have a descent in 
position $m$ if and only if  $m$ is divisible by  $k$.
The classical Euler numbers are the special case when $k=2$.
In this paper, we  study  divisibility properties of a $q$-analog of
$E_{n|k}$.  In particular, we generalize two theorems of
Andrews and Gessel~[3] about factors of the $q$-tangent numbers.
         \end{abstract}
\pagestyle{plain}

\section{Introduction}                                    \label{i}

Let $\fS_n$ denote the {\it symmetric group} of all permutations of
the set $\bn=\{1,2,\ldots,n\}$.
The $n$th {\it Euler number}, $E_n$, can be defined as the number of
permutations $\pi=a_1a_2\ldots a_n$ in $\fS_n$
that alternate, i.e.,
$$
a_1<a_2>a_3<\ldots
$$
These numbers have a long and venerable history going back at least to
Andr\'e~[1,2].  Comtet's book~[4, p.\ 48] lists some of the 
classical properties of the $E_n$.  In particular, the Euler numbers
have exponential generating function
$$
\sum_{n\ge0} E_n x^n/n!=\tan x + \sec x.
$$
For this reason the $E_{2n+1}$ are called {\it tangent numbers} and
the $E_{2n}$ {\it secant numbers}.

The {\it descent set} of any
$\pi=a_1a_2\ldots a_n$ is the set of indices 
$$
\Des(\pi)=\{i\ :\ 1\le i<n \mbox{ and } a_i>a_{i+1}\}.
$$
Furthermore, the  {\it generalized Euler}, $E_{n|k}$, counts the
number of $\pi\in\fS_n$ such that $\Des(\pi)=\{k,2k,3k,\ldots\}$.
We will also use $E_{\bn|k}$ to denote the set
of all such permutations.  Clearly, $E_{n|2}=E_n$.  
As an example, we have
$$
E_{{\bf 5}|3}=
\{12435,\ 13425,\ 23415,\ 12534,\ 13524,\ 14523,\ 23514,\ 24513,\ 34512\}.
$$

We will be
concerned with a certain $q$-analog of the generalized Euler numbers
defined as follows.  An {\it inversion} of $\pi=a_1a_2\ldots a_n$ is
an out-of-order pair, namely $(a_i,a_j)$ with $i<j$ and $a_i>a_j$.  We
let $\inv\pi$ denote the number of inversions of $\pi$.  Following
Stanley~[7, pp.\ 147--9], define
\beq                                                    \label{Eq}
      E_{n|k}(q) = \sum_{\pi} q^{\inv(\pi)}.
\eeq
Continuing our example from the previous paragraph,
\bea
E_{5|3}(q)&=&q+q^2+q^3+q^2+q^3+q^4+q^4+q^5+q^6\\
	&=&q+2q^2+2q^3+2q^4+q^5+q^6.
\eea

It is well-known that the tangent numbers are divisible by high powers
of 2.  In~[3], Andrews and Gessel show that 
both $(1+q)(1+q^2)\cdots(1+q^n)$ and $(1+q)^n$ divide
the $q$-tangent number $E_{2n+1}(q)$.  Our main theorem is a
generalization of this result to $E_{n|k}(q)$.  Let 
$$
[k]=[k]_q=1+q+q^2+\cdots+q^{k-1}.
$$
So $[k]_{q^i}=1+q^i+q^{2i}+\cdots+q^{(k-1)i}$.
\bth
Let $k$  be prime and $ 1 \leq i \leq k-1$.  Then
\ben
\item $[k] [k]_{q^2} [k]_{q^3} \cdots [k]_{q^n}\ |\ E_{nk+i|k}(q)$;
\item $[k]^n\ |\ E_{nk+i|k}(q)$.
\een
\eth

The next section is devoted to proving a recursion and two lemmas that
we will need for the proof of the
previous theorem.  Section 3 is devoted to the demonstration of the
theorem itself.  Finally, we close with a section of comments and open
questions.

\section{Lemmas}					\label{l}

It will be useful to have a recursion relation for the
$E_{n|k}(q)$.  To state it we will need the {\it Gaussian polynomials}
or {\it $q$-binomial coefficients}
$$
\gau{n}{k}=\frac{[n]!}{[k]![n-k]!}
$$
where $[n]!=[n][n-1]\cdots[1]$.   We assume $\gau{n}{k}=0$ if $k>n$.
It is well-known that these polynomials can be written as
$$
\gau{n}{k}=\sum_{\pi\in\fS_{n,k}} q^{\inv\pi}
$$
where $\fS_{n,k}$ is the set of all permutations of $k$ zeros and
$n-k$ ones.  Finally, let $\chi(P)$ be the {\it characteristic
function} which is $1$ if the statement $P$ is true and $0$ if it is false.


\bpr                                                         \label{qr}
For $n\ge0$ the $E_{n|k}(q)$ satisfies the  following  recursion
$$
E_{(n+1)|k}(q)
   =\sum_{m=1}^{\fl{n/k}}
\gau{n}{mk-1} q^{n-mk+1}E_{(mk-1)|k}(q)E_{(n-mk+1)|k}(q)
+\chi(k\not|\ n)E_{n|k}(q)
$$
with boundary condition $E_{0,k}(q)=1$.
\epr
\pf\
The initial condition is trivial.  For the recurrence relation,
consider all the indices $i$ where one could have $a_i=n+1$ in some
$\pi=a_1a_2\ldots a_{n+1}\in\fS_{n+1}$.  Clearly, $i=n+1$ can occur iff
$k\not|\ n$ and in this case all $\pi$ ending in $n+1$ contribute
$E_{n|k}(q)$ to the sum~(\ref{Eq}) for $E_{n+1|k}$.  

The other possible positions for $n+1$ are at $i=mk$ for some $m$,
$1\le m\le\fl{n/k}$.  In this case, the inversions caused by $n+1$ are
accounted for by $q^{n-km+1}$.  The inversions of the elements
$a_1\ldots a_{mk-1}$ and $a_{mk+1}\ldots a_{n+1}$ among themselves are
counted by $E_{(mk-1)|k}(q)$ and $E_{(n-mk+1)|k}(q)$, respectively.
Finally, $\gau{n}{mk-1}$ takes care of inversions between these two
sets of integers. \hfill\Qed

Given the form of the recursion in the previous proposition, it should
come as no surprise that we will
need two lemmas about divisibility properties of $q$-binomial coefficients .
\ble 								\label{qc2}
Let $k$ be a prime  and suppose $ 0 \leq i \leq k-2$.
Then for any non-negative 
integers n and m, we have the divisibility relation
$$
[k]\ |\ \gau{nk+i}{mk-1}.
$$
\ele
\pf\  For all non-negative integers $n$ and $m$, we have
\beq						\label{eq}
\gau{nk+i}{mk-1}=
\frac{[nk+i][nk+i-1]\cdots[nk-mk+i+2]}{[mk-1]!}.
\eeq
Since $k$ is prime, all roots of $[k]$ are primitive $k$th roots of unity.
If $\ze$ is such a root then $\ze$ is a root of $[l]$ iff $k|l$ and in
that case it has multiplicity one.
Thus the multiplicity of $\ze$ as a root of the numerator of~(\ref{eq})  is
$n-(n-m)=m$ while in the denominator it is
$m-1$.  Thus $[k]$ divides the $q$-binomial coefficient as claimed.\hfill \qed

\ble                                                        \label{qc1} 
Let  $k$ be prime and suppose $0\leq i \leq k-2$.
Then for any non-negative integers n and m, the expression
$$
\gau{nk+i}{mk-1}
\frac{[k]_{q^{m-1}}[k]_{q^{m-2}}\cdots[k]}
{[k]_{q^n}[k]_{q^{n-1}}\cdots[k]_{q^{n-m+1}}} 
$$
is a polynomial in $q$. 
\ele
\pf\ As with the previous lemma,
we need only show that each root of unity which is a zero of the
denominator appears with at least as large multiplicity in the numerator.
We write the Gaussian polynomial as
     \beq 						\label{quot1}
	\gau{nk+i}{mk-1}=
      \frac{(1-q^{nk+i} ) (1-q^{nk+i-1} )\cdots(1-q^{nk-mk+i+2} )} 
     {(1-q^{mk-1} ) (1-q^{mk-2} )\cdots(1-q^2 ) (1-q)}=
	P_1 P_2
     \eeq
where
$$
P_1=\frac{(1-q^{nk})(1-q^{(n-1)k})\cdots(1-q^{(n-m+1)k})}
       {(1-q^{(m-1)k})(1-q^{(m-2)k})\cdots(1-q^{k})}
$$
and $P_2$ contains all the other factors of $1-q^j$.
Substituting $[k]_{q^j} = (1-q^{jk} )/(1-q^j)$ into the expression in the
statement of the lemma and doing some cancelation shows that
\beq							\label{quot2}
\gau{nk+i}{mk-1}
\frac{[k]_{q^{m-1}}[k]_{q^{m-2}}\cdots[k]}
{[k]_{q^n}[k]_{q^{n-1}}\cdots[k]_{q^{n-m+1}}} 
=P_3 P_2
\eeq
where
$$
P_3=\frac{(1-q^n)(1-q^{n-1})\cdots(1-q^{n-m+1})}
{(1-q^{m-1})(1-q^{m-2})\cdots(1-q)}.
$$

Let $\ze$ be a primitive $l$th root of unity.  Then since $k$ is
prime, either $\gcd(k,l)=1$ or $k|l$.  In the former case, the
multiplicities of $\ze$ in the denominators of $P_1P_2$ and $P_3P_2$
are equal.  Since the same is true of the numerators, $P_3P_2$ has no
pole at $\ze$ since $P_1P_2$ doesn't.  If $k|l$ then $\ze$ is neither
a root nor a pole of $P_2$.  Also $P_3=(1-q^n)\gau{n-1}{m-1}$ is a
polynomial in $q$ and so does not have $\ze$ as a pole.  Thus no
primitive root of unity is a pole of $P_2P_3$ forcing it to be a
polynomial. \hfill\qed

\section{Divisibility of the  generalized $q$-Euler numbers}

Now we are in the position to prove our main results.

\bth                                                \label{qd1} 
Let k be prime and $1 \leq i \leq k-1$, then
$E_{(nk+i)|k}(q)$ is divisible by $[k]^n$.
\eth
\pf\ We will induct on $n$.  For $n=0$, the result is trivial. 
Suppose the result  is true up to but not including $n$.
First consider $i=1$.  According to Proposition~\ref{qr}
\beq						\label{i=1}
E_{(nk+1)|k}(q)
   =\sum_{m=1}^{n}
q^{nk-mk+1}\gau{nk}{mk-1} E_{(mk-1)|k}(q)E_{(nk-mk+1)|k}(q).
\eeq
By induction 
$[k]^{m-1}$ and $[k]^{n-m}$ divide $E_{(mk-1)|k}(q)$ and
$E_{(nk-mk+1)|k}(q)$, respectively.  But
by Lemma \ref{qc2}, $[k]$ is a factor of the corresponding
$q$-binomial coefficient in~(\ref{i=1}).  So $E_{(nk+1)|k}(q)$ is
divisible by $[k]^{m-1+n-m+1}=[k]^n$ as desired.
The case when $2\le i\le k-1$ is similar, with the extra term from the 
recursion in Proposition~\ref{qr} being taken care of by the case for
$i-1$.\hfill \qed

\bth							\label{qd2}
Let $k$  be prime and $ 1 \leq i \leq k-1$, then
$E_{(nk+i)|k}(q)$ is divisible by
$[k] [k]_{q^2} [k]_{q^3}\cdots[k]_{q^n}$.
\eth
\pf\  As before, we will induct on $n$ with $n=0$ being trivial. 
Suppose the result is true up to but not including $n$. 
When $i=1$ we have~(\ref{i=1}) again and examine each of its terms.
By the induction hypothesis
   $$E_{(mk-1)|k}(q) = [k][k]_{q^2}[k]_{q^3}\cdots[k]_{q^{m-1}} Q_1$$
and
   $$E_{(nk-mk+1)|k}(q)=[k] [k]_{q^2} [k]_{q^3}\cdots[k]_{q^{n-m}} Q_2$$
where $Q_1$ and $Q_2$ are  polynomials in $q$.   Then
$$
\barr{l}
{\displaystyle \gau{nk}{mk-1} E_{(mk-1)|k}(q)E_{(nk-mk+1)|k}(q)}\\[10pt]
{\displaystyle\rule{20pt}{0pt}\gau{nk}{mk-1}
\frac{[k]_{q^{m-1}}[k]_{q^{m-2}}\cdots[k]}
{[k]_{q^n}[k]_{q^{n-1}}\cdots[k]_{q^{n-m+1}}} 
          [k] [k]_{q^2} [k]_{q^3} \cdots[k]_{q^n} Q_1 Q_2.}
\earr
$$
By Lemma~\ref{qc1}, the $q$-binomial coefficient times the fraction 
is a polynomial in $q$.  So $[k] [k]_{q^2} [k]_{q^3} \cdots[k]_{q^n}$ 
is a factor of every term in~(\ref{i=1}) and thus of $E_{(nk+1)|k}(q)$.
The case $2\le i\le k-1$ is handled as in the proof of
Theorem~\ref{qd1}.\hfill\qed 

\section{Comments and open questions}

By setting $q=1$ in either Theorem~\ref{qd1} or Theorem~\ref{qd2} we get the
following corollary.
\bco						
Let $k$  be prime and $ 1 \leq i \leq k-1$, then
$E_{(nk+i)|k}$ is divisible by
$k^n$.
\eco
It is well known that for $k=2$ (the tangent numbers) 
\beq							\label{tan}
2^{2n}\ |\ (n+1)E_{2n+1}
\eeq 
and that the corresponding quotient, called a {\it
Genocchi number}, is odd.  Thus it is not surprising that better
divisibility results can be obtained when $q=1$ for general primes $k$.  In
particular, Gessel and Viennot~[6] have shown that
\beq							\label{tank}
k^{\cef{nk-j}{k-1}}\ |\ {nk\choose j}E_{(nk-j)|k}.
\eeq
Note that this reduces to~(\ref{tan}) when $k=2$ and $j=1$.  This raises a
couple of questions.  Is it true that the associated quotient
in~(\ref{tank}) is relatively prime to $k$?  Can these results be
extended to the case of 
arbitrary $q$?  With regards to the second query, the reader should
consult Foata's article~[5] which provides some answers in
the case $k=2$.

\hs{.5in}

\noindent{\bf\Large References}

\ben
\item Andr\'e, D.: D\'eveloppements de $\sec x$ et de
$\tang x$, {\it C. R. Acad.\ Sci.\ Paris}, {\bf 88}, 965--967 (1879).
\item Andr\'e, D.: Sur les permutations altern\'ees,
{\it J. Math.\ Pures Appl.,} {\bf 7}, 167--184 (1881).
\item Andrews, G., Gessel, I.: Divisibility
properties of the $q$-tangent numbers, {\it \pams}, {\bf 68}, 380--384 (1978).
\item Comtet, L., {\it Advanced Combinatorics,} D. Reidel
Pub.\ Co., 1974.
\item Foata, D.: Further divisibility properties of
the $q$-tangent numbers, {\it \pams}, {\bf 81}, 143--148 (1981).
\item Gessel, I., Viennot, G.: Generalized tangent and
Genocchi numbers, preprint.
\item Stanley, R.P., {\it Enumerative Combinatorics,
Volume I,} Wadsworth and Brooks/Cole, 1986.
\een
\end{document}